\documentclass[runningheads]{lncse}

\usepackage{amsmath,amssymb}
\usepackage{amsfonts}
\usepackage{graphicx}
\usepackage[usenames, dvipsnames]{color}

\graphicspath{ {./eps/} }

\begin{document}

\title{Model Order Reduction for Pattern Formation in
FitzHugh-Nagumo Equation}

\titlerunning{Model Order Reduction for Pattern Formation}

\author{ B\"{u}lent Karas\"{o}zen\inst{1}  \and Murat Uzunca\inst{2} \and Tu\u{g}ba K\"u\c{c}\"ukseyhan\inst{1}  }

\authorrunning{Karas\"{o}zen et al.}

\institute{
Department of Mathematics and Institute of Applied Mathematics, Middle East Technical University, 06800 Ankara, Turkey, {\tt bulent@metu.edu.tr}, {\tt guney.tugba@metu.edu.tr} \and Institute of Applied Mathematics, Middle East Technical University, 06800 Ankara, Turkey, {\tt uzunca@gmail.com}
}

\maketitle

\begin{abstract}
We developed a reduced order model (ROM) using the proper orthogonal decomposition (POD)  to compute efficiently the labyrinth and spot like patterns of the FitzHugh-Nagumo (FNH)  equation. The FHN equation is discretized in space by the discontinuous Galerkin (dG) method  and in time by the backward Euler method. Applying POD-DEIM (discrete empirical interpolation method) to the full order model (FOM) for different values of the parameter in the bistable nonlinearity,  we show that using few POD and DEIM modes, the patterns can be computed accurately. Due to the local nature of the dG discretization, the POD-DEIM requires less number of connected nodes than continuous finite element for the nonlinear terms, which leads to a significant reduction of the computational cost for dG POD-DEIM.
\end{abstract}

\section{Introduction}
There has been significant development in the efficient implementation and analysis of the model order reduction (MOR) techniques for parametrized partial differential equations (PDEs) \cite{hesthaven15}. Even though the POD is a very successful MOR technique  for linear problems, for nonlinear problems the computational complexity of the evaluation of the nonlinear reduced model still depends on the dimension of the FOM. Several methods are developed to reduce the computational cost so that the nonlinear function evaluations are independent of the dimension of the FOM  and the computational complexity is proportional to the dimension of ROM.  The discrete empirical interpolation method (DEIM) \cite{chaturantabut10} which is the modified version of the  empirical interpolation method (EIM) \cite{barrault04} are the most frequently used ones. The DEIM was originally developed for nonlinear functions which depend component-wise on single variables, arising from the finite difference discretization of nonlinear PDEs \cite{chaturantabut10}. When the nonlinear functions are discretized by finite elements,  the discretized nonlinear functions depend on the mesh and on the polynomial degree of the finite elements. Therefore the efficiency of the POD-DEIM can be degraded. In \cite{antil14} the DEIM was  applied at different stages of the finite element assembly process. Using the unassembled finite elements, where each DEIM point is related to a single element, the number of nonlinear function calls during the online computation is reduced, but the size of the nonlinear snapshots are enlarged, which increases the offline computational cost \cite{antil14}. In  this paper we consider the dG discretization for time dependent parametrized semi-linear PDEs. Due to the local nature of the dG approximation, each component of the discretized nonlinear vectors depend only on few elements in the local mesh, whereas the continuous FEMs discretized nonlinear vectors depend on multiple components of the finite element solutions. Therefore the number of POD-DEIM function evaluations for dG approximation  is comparable with the finite difference discretization.

In this paper we consider the parametrized FHN equation \cite{lago14} with fast diffusing inhibitor $D_v > D_u$
\begin{subequations}\label{fhn}
\begin{align}
\frac{\partial u}{\partial t} & = D_u\Delta u  - \alpha (v-u) - f(u;\mu ), \\
\frac{\partial v}{\partial t} & = D_v\Delta u  - \beta (v-u),
\end{align}
\end{subequations}
on a space-time cylinder $\Omega\in \mathbb{R}^2\times (0,T]$ with homogeneous (zero-flux) Neumann boundary conditions. The variables $u(x,t;\mu)$ and $v(x,t;\mu)$ stand for the activator and inhibitor, respectively. The term $f(u;\mu ) = (u-\mu)(u^2-1)$ represents the bistable nonlinearity for the parameter $\mu$. We investigate the formation of labyrinth and spot like patterns for different values of the parameter $\mu$ as in \cite{lago14}, where  $\alpha$, $\beta$ and diffusion coefficients  $D_u$, $D_v$ are fixed.

The paper is organized as follows. In Section 2 the discretization of FHN equation \eqref{fhn} in space by symmetric interior penalty discontinuous Galerkin method  (SIPG) is given. In Section 3 the ROM based on the POD is formulated. In  Section 4 we describe the SIPG discretized version of the DEIM for the bistable nonlinearity. Numerical results for pattern formations for different values of the parameters demonstrate the good performance of dG ROMs. The paper ends with some conclusions and outlook for the future work.

\section{Full Order Model}

  The FHN equation \eqref{fhn} is discretized in space using SIPG method \cite{riviere08dgm}. Let $\varepsilon_h$ be the disjoint partition of the domain $\Omega\subset\mathbb{R}^2$ with elements (triangles) $\{E_i\}_{i=1}^{N_{el}}\in\varepsilon_h$, where $N_{el}$ is the number of elements in the partition. The discrete solution and test function spaces on $\varepsilon_h$ are given by
\begin{equation*}
D_q=D_q(\varepsilon_h):=\{w \in L^2(\Omega): w_E \in \mathbb{P}_q(E),\; \forall E \in \varepsilon_h\},
\end{equation*}
where $\mathbb{P}_q(E)$ is the space of polynomials of degree at most $q$ on $E\in \varepsilon_h$, and the functions $w\in D_q$ are discontinuous along the inter-element boundaries. Multiplying  \eqref{fhn} by arbitrary test functions $w_1, w_2\in D_q$ and integrating by using Green's theorem over each mesh element, we obtain the semi-discrete variational equations
\begin{equation} \label{dgvar}
\begin{aligned}
\left(\frac{\partial u_h}{\partial t}, w_1\right) + a_h(D_u;u_h, w_1) + \alpha( v_h-u_h,w_1) + (f(u_h;\mu ), w_1) &= 0, \\
\left(\frac{\partial v_h}{\partial t}, w_2\right) + a_h(D_v;u_h, w_2) + \beta( v_h-u_h,w_2) &= 0,
\end{aligned}
\end{equation}
where $a_h(D_u ;u ,w_1)$ and  $a_h(D_v ;u ,w_2)$ refer to the dG bilinear forms. We refer to  \cite{riviere08dgm} for the details of dG discretization.

Let $N:=N_{loc}\times {N_{el}}$  denotes the dG degrees of freedom (DoFs), where $N_{loc}$ is the local dimension on each element depending on the polynomial degree $q$. Then, for any $t\in(0,T]$, the dG solutions of \eqref{dgvar} are of the form
\begin{equation}\label{fhn_comb}
u_h(t,x)=\sum_{i=1}^{N}u_i(t)\phi_i(x)=\phi\mathbf{u} \; , \quad v_h(t,x)=\sum_{i=1}^{N}v_i(t)\phi_i(x)=\phi\mathbf{v},
\end{equation}
where ${\mathbf u}(t):=(u_1(t),\ldots , u_{N}(t))^T$ and ${\mathbf v}(t):=(v_1(t),\ldots , v_{N}(t))^T$ are the vectors of time dependent unknown coefficients of $u_h$ and $v_h$, respectively, and $\phi (x):=[\phi_1(x)\; \ldots \; \phi_{N}(x)]$ is the matrix of the basis functions. Plugging \eqref{fhn_comb} into the equations \eqref{dgvar} and choosing $w_1=w_2=\phi_i$, $i=1,\cdots,N$,  we obtain the FOM of \eqref{fhn} as the following system of ordinary differential equations (ODEs)
\begin{equation}\label{fhn_ode}
\begin{aligned}
M {\mathbf u}_t + S_u {\mathbf u} + \alpha M({\mathbf v}-{\mathbf u}) + F({\mathbf u};\mu ) &= 0, \\
M {\mathbf v}_t + S_v {\mathbf v} + \beta M({\mathbf v}-{\mathbf u}) &= 0,
\end{aligned}
\end{equation}
where $S_u,S_v\in\mathbb{R}^{N\times N}$ are the stiffness matrices, $M\in \mathbb{R}^{N\times N}$ is the mass matrix and $F({\mathbf u};\mu )\in \mathbb{R}^N$ is the nonlinear vector depending on the parameter $\mu$. The FOM \eqref{fhn_ode} is solved in time by the backward Euler method.

\section{Reduced Order Model}

For an arbitrary parameter $\bar{\mu}$, the $k$-th order approximate ROM solutions $\tilde{u}_{h,k}:=\tilde{u}_{h,k}(t,x;\bar{\mu})$ and $\tilde{v}_{h,k}:=\tilde{v}_{h,k}(t,x;\bar{\mu})$ have the form
\begin{equation}\label{fhn_combrom}
\tilde{u}_{h,k}= \sum_{i=1}^{k} \tilde{u}_i(t) \psi_{u,i}(x) \; , \quad \tilde{v}_{h,k}= \sum_{i=1}^{k} \tilde{v}_i(t) \psi_{v,i}(x),
\end{equation}
where $\mathbf{\tilde{u}}(t):=(\tilde{u}_1(t),\ldots ,\tilde{u}_k(t))^T$ and $\mathbf{\tilde{v}}(t):=(\tilde{v}_1(t),\ldots ,\tilde{v}_k(t))^T$ are the coefficient vectors of the ROM solutions. For a set $\{\mu_1,\ldots , \mu_{n_s}\}$ of parameter samples, the POD  reduced basis functions $\{\psi_{u,i}\}$ and $\{\psi_{v,i}\}$ are computed  as  solutions of the minimization problem
\begin{equation}\label{min}
\begin{aligned}
\min_{\psi_{w,1},\ldots ,\psi_{w,k}} \frac{1}{n_s}\sum_{m=1}^{n_s}\frac{1}{J}\sum_{j=1}^J \left\| w^{m,j} - \sum_{i=1}^k (w^{m,j},\psi_{w,i})_{L^2(\varepsilon_h)}\psi_{w,i}\right\|_{L^2(\varepsilon_h)}^2 \\
\text{subject to } (\psi_{w,i},\psi_{w,j})_{L^2(\varepsilon_h)} = \Psi_{w,\cdot ,i}^TM\Psi_{w,\cdot ,j}=\delta_{ij} \; , \; 1\leq i,j\leq k,
\end{aligned}
\end{equation}
for $w^{m}\in\{u(t,x;\mu_m),v(t,x;\mu_m)\}$, $m=1,\ldots , n_s$. In \eqref{min}, the terms $\{ w^{m,j}\}_{j=1}^J$ denote the approximate solutions $w^{m,j}\approx w^{m}(t_j)$ at the time  $t_j$ for a fixed parameter $\mu_m$ and $\delta_{ij}$ is the Kronecker delta. We note that $w^{m,j}$ in \eqref{min} stand for the solutions but not for the coefficients of the unknown solutions as in the continuous finite elements. For the dG discretizations we use  modal basis functions where the coefficients do not coincide with the solution values. Therefore, instead of the Euclidean norm, we use in \eqref{min} the weighted inner product and the corresponding norm with the symmetric positive definite mass matrix $M$, leading to $M$-orthogonal reduced basis functions \cite{hinze08}.

In practice, instead of the minimization problem \eqref{min}, the equivalent eigenvalue problem is solved \cite{hinze08}.
\begin{equation}\label{eg2}
\widehat{\mathcal{U}}\widehat{\mathcal{U}}^T\widehat{\Psi}_{u,\cdot ,i}=\sigma_{u,i}^2\widehat{\Psi}_{u,\cdot ,i} \; , \qquad \widehat{\mathcal{V}}\widehat{\mathcal{V}}^T\widehat{\Psi}_{v,\cdot ,i}=\sigma_{v,i}^2\widehat{\Psi}_{v,\cdot ,i} \; , \quad i=1,2,\ldots ,k,
\end{equation}
where $\widehat{\mathcal{U}}=R\mathcal{U}$, $\widehat{\mathcal{V}}=R\mathcal{V}$, $\widehat{\Psi}_{\cdot,\cdot ,i}=R\Psi_{\cdot,\cdot ,i}$, $R^T$ is the Cholesky factor of the mass matrix $M$, and $\mathcal{U}=[\mathbf{u}^{1,1},\ldots ,\mathbf{u}^{n_s,J}]$ and $\mathcal{V}=[\mathbf{v}^{1,1},\ldots ,\mathbf{v}^{n_s,J}]$ in $\mathbb{R}^{N\times (n_s\times J)}$ are the snapshot matrices. The vectors $\Psi_{u,\cdot ,i}$ and $\Psi_{v,\cdot ,i}$ denote the coefficient vectors of the reduced basis functions $\psi_{u,i}$ and $\psi_{v,i}$, respectively. The solutions $\widehat{\Psi}_{\cdot,\cdot ,i}$  of \eqref{eg2} are obtained as the first $k$ left singular vectors in the generalized singular value decomposition (SVD)  of $\widehat{\mathcal{U}}$ and $\widehat{\mathcal{V}}$, respectively \cite{hinze08}. Combining the FOM  \eqref{fhn_comb} and ROM \eqref{fhn_combrom} solutions and using  the fact that the reduced basis functions $\psi_{u,i}$ and $\psi_{v,i}$ belong to the dG space $D_q$, we obtain the relations between the coefficient vectors $\mathbf{u}$, $\mathbf{v}$ of the FOM solutions and the coefficient vectors $\mathbf{\tilde{u}}$, $\mathbf{\tilde{v}}$ of the ROM solutions
\begin{align}\label{fulltorom}
\mathbf{u} = \Psi_u\mathbf{\tilde{u}} \; , \quad \mathbf{v} = \Psi_v\mathbf{\tilde{v}}.
\end{align}
Substituting   \eqref{fulltorom} into  \eqref{fhn_ode} and projecting onto the reduced spaces spanned by $\{\psi_{u,1},\ldots ,\psi_{u,k}\}$  and $\{\psi_{v,1},\ldots ,\psi_{v,k}\}$, respectively, we obtain  the $k$-dimensional ROM
\begin{equation}\label{fhn_rom}
\begin{aligned}
\mathbf{\tilde{u}}_t  +  \tilde{S}_u\mathbf{\tilde{u}} + \alpha \tilde{M}_u\mathbf{\tilde{v}} - \alpha\mathbf{\tilde{u}} + \Psi_u^TF(\Psi_u\mathbf{\tilde{u}};\bar{\mu})  &= 0,\\
\mathbf{\tilde{v}}_t  +  \tilde{S}_v\mathbf{\tilde{v}} + \beta\mathbf{\tilde{v}} - \beta\tilde{M}_v\mathbf{\tilde{u}}  &= 0
\end{aligned}
\end{equation}
with the reduced matrices
$$
\tilde{S}_u=\Psi_u^TS_u\Psi_u \; , \quad \tilde{S}_v=\Psi_v^TS_v\Psi_v \; , \quad \tilde{M}_u=\Psi_u^TM\Psi_v\; , \quad \tilde{M}_v=\Psi_v^TM\Psi_u.
$$
The system \eqref{fhn_rom} is solved by backward Euler method, as well.

\section{Discrete Empirical Interpolation Method (DEIM)}

Although the dimension of the reduced system \eqref{fhn_rom} is small, $k\ll N$, the computation of the nonlinear term $N(\mathbf{\tilde{u}}):=\Psi_u^TF(\Psi_u\mathbf{\tilde{u}};\bar{\mu})$ still depends on the dimension $N$ of the full system. We apply DEIM \cite{chaturantabut10} to reduce the computational cost, where the nonlinear function is approximated as $F(\Psi_u\mathbf{\tilde{u}};\bar{\mu}) \approx Ws(t)$, from a subspace $W=[W_1,\ldots ,W_n]\in\mathbb{R}^{N\times n}$, where each member $W_i$ is called the DEIM basis functions, $i=1,2,\ldots ,n$ ($n\ll N$). The DEIM basis functions $W_i$ are computed through the SVD of the nonlinear snapshot matrix  $\mathcal{F}:=[F^{1,1},\ldots , F^{n_s,J}]\in\mathbb{R}^{N\times (n_s\times J)}$, where $F^{m,i}:=F(\Psi_u\mathbf{\tilde{u}}(t_i);\mu_m)$ are the nonlinear vectors at the time instance $t_i$, obtained in the online computation for the parameters $\mu_m$, $m=1,\ldots ,n_s$. Because the system $Ws(t)$ is overdetermined, the  projection matrix $P$ is introduced which is computed by the greedy DEIM algorithm \cite{chaturantabut10}. Then, we  use the approximation $N(\mathbf{\tilde{u}}) \approx \tilde{N}(\mathbf{\tilde{u}}) = Q\tilde{F}^{\bar{\mu}}$ where the matrix $Q=\Psi_u^TW(P^TW)^{-1}\in\mathbb{R}^{k\times n}$ is precomputable and $\tilde{F}^{\bar{\mu}}=P^T F(\Psi_u\mathbf{\tilde{u}};\bar{\mu})\in\mathbb{R}^{n}$ is the $n$-dimensional non-linear vector which can be computed in an efficient way. In addition, the DEIM approximation satisfies the a priori  error bound 
$$
\|F^{\bar{\mu}}-W(P^TW)^{-1}\tilde{F}^{\bar{\mu}}\|_2\leq \|(P^TW)^{-1}\|_2\|(I-WW^T)F^{\bar{\mu}}\|_2,
$$
where the term $\|(P^TW)^{-1}\|_2$ is of moderate size of order 100 or less \cite{antil14}.

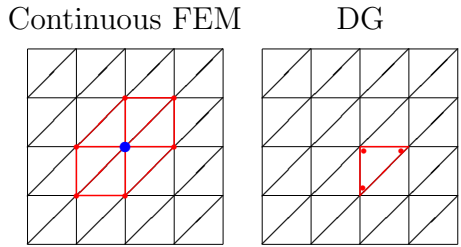
\begin{figure}[htb]
	\centering
	\setlength{\unitlength}{0.26mm}
		\begin{picture}(220, 120)
		
				\linethickness{0.1mm}
				\put(-10,110){\large{Continuous FEM}}
				\put(0,0){\line(1,0){100}}
				\put(0,0){\line(0,1){100}}
				\put(0,100){\line(1,0){100}}
				\put(100,0){\line(0,1){100}}
								
				\put(0,25){\line(1,0){100}}
				\put(0,50){\line(1,0){100}}
				\put(0,75){\line(1,0){100}}
				\put(25,0){\line(0,1){100}}
				\put(50,0){\line(0,1){100}}
				\put(75,0){\line(0,1){100}}
				
				\put(0,0){\line(1,1){100}}
				\put(0,25){\line(1,1){75}}
				\put(0,50){\line(1,1){50}}
				\put(0,75){\line(1,1){25}}
				\put(25,0){\line(1,1){75}}
				\put(50,0){\line(1,1){50}}
				\put(75,0){\line(1,1){25}}
				
				\color{red}
				\linethickness{0.2mm}
				\put(50,25){\line(1,1){25}}				
				\put(50,25){\line(0,1){25}}
				\put(50,50){\line(1,0){25}}			
				\put(25,25){\line(1,0){25}}	
				\put(25,25){\line(1,1){25}}	
				\put(25,25){\line(0,1){25}}	
				\put(25,50){\line(1,0){25}}	
				\put(25,50){\line(1,1){25}}	
				\put(50,50){\line(0,1){25}}		
				\put(50,50){\line(1,1){25}}		
				\put(50,75){\line(1,0){25}}		
				\put(75,50){\line(0,1){25}}	
					
				\put(25,25){\circle*{2.5}}
				\put(50,25){\circle*{2.5}}
				\put(25,50){\circle*{2.5}}
				\put(75,50){\circle*{2.5}}
				\put(50,75){\circle*{2.5}}
				\put(75,75){\circle*{2.5}}
				\color{blue}
				\put(50,50){\circle*{5}}

				\color{black}
				\linethickness{0.1mm}
				\put(158,110){\large{DG}}
				\put(120,0){\line(1,0){100}}
				\put(120,0){\line(0,1){100}}
				\put(120,100){\line(1,0){100}}
				\put(220,0){\line(0,1){100}}
				
				\put(120,25){\line(1,0){100}}
				\put(120,50){\line(1,0){100}}
				\put(120,75){\line(1,0){100}}
				\put(145,0){\line(0,1){100}}
				\put(170,0){\line(0,1){100}}
				\put(195,0){\line(0,1){100}}

				\put(120,0){\line(1,1){100}}
				\put(120,25){\line(1,1){75}}
				\put(120,50){\line(1,1){50}}
				\put(120,75){\line(1,1){25}}
				\put(145,0){\line(1,1){75}}
				\put(170,0){\line(1,1){50}}
				\put(195,0){\line(1,1){25}}
				
				\color{red}
				\linethickness{0.2mm}
				\put(170,25){\line(1,1){25}}				
				\put(170,25){\line(0,1){25}}
				\put(170,50){\line(1,0){25}}	
					
			  \put(171.5,29){\circle*{2.5}}
				\put(172,48){\circle*{2.5}}
				\put(191,48){\circle*{2.5}}
				
				\color{black}
				\linethickness{0.1mm}
		\end{picture}
		\caption{Connectivity of degrees of freedoms for linear basis functions.\label{connect}}		
	\end{figure}

In dG discretization, the integrals are computed on a single triangular element, whereas  for  continuous finite element discretizations with linear polynomials all the interior degrees of freedoms are shared by  6 triangular elements, see Fig.~\ref{connect}. The unassembled finite element approach is used in \cite{antil14}, so that each DEIM point is related to one element, which reduces the online computational cost,  but  increases the number of snapshots and therefore the cost of the offline computation. Due to its local nature, the dG discretization is automatically in the unassembled form and it does not require computation of additional snapshots.

\vspace*{- 0.5 cm}
\vspace*{- 0.5 cm}
As we use for the time integration the implicit backward Euler method, on each time step the nonlinear equations have to be solved by Newton's method. Therefore the computational cost of the Jacobian by DEIM of the  reduced model has to be  taken into account. Because the support of dG basis functions have only a single element, the Jacobian matrix of the FOMs appears in block diagonal form unlike the continuous FEMs where the Jacobian matrix contains overlapping blocks. The Jacobian matrix arising from POD and POD-DEIM are of the form
$$
\frac{\partial}{\partial \mathbf{\tilde{u}}} N(\mathbf{\tilde{u}}) = \Psi_u^T J_F^{\bar{\mu}} \Psi_u \; , \quad \frac{\partial}{\partial \mathbf{\tilde{u}}} \tilde{N}(\mathbf{\tilde{u}}) = Q (P^TJ_F^{\bar{\mu}}) \Psi_u,
$$
where $(P^TJ_F^{\bar{\mu}})\in\mathbb{R}^{n\times N}$ is the matrix including only $n\ll N$ rows of the Jacobian $J_F^{\bar{\mu}}$, and in each row of the Jacobian there are only $N_{loc}$ nonzero terms because of the local structure of the dG. Hence, only $n\times N_{loc}$ entries are needed to compute $P^TJ_F^{\bar{\mu}}$, whereas without DEIM,  $N_{el}\times N_{loc}^2$ entries are  required for computation of the Jacobian $J_F^{\bar{\mu}}$ of the FOM.


\section{Numerical Results }
We consider FHN equation \eqref{fhn} for $(x,t)\in [-10,10]^2\times [0,1000]$ with  random  initial conditions uniformly distributed between $-1$ and $1$. The  other parameters  $D_u=0.04$, $D_v=1$, $\alpha =0.3$, $\beta =1$ are fixed as in   \cite{lago14}. We use linear dG polynomials ($N_{loc}=3$), and as the discrete mesh, we form the partition of $[-10,10]^2$, by 5 times uniform refinement,  with 2048 triangular elements leading to 6144  DoFs. Snapshots are taken in the time interval $[0,1000]$ with the time step $\Delta t=0.5$. For POD/POD-DEIM basis construction, we use the parameter samples $\mu\in\{ -0.04,-0.02,0,0.02,0.04\}$, $n_s=5$. The reduced systems are solved for the set $\{ -0.03,-0.01,0.01,0.03\}$ of parameter values of $\mu$, which are not contained in the set of sample parameters. The average number of Newton iterations was 1 for the computation of the FOMs and ROMs on each time step.

\begin{figure}[htb!]
\centering
\includegraphics[width=0.35\textwidth]{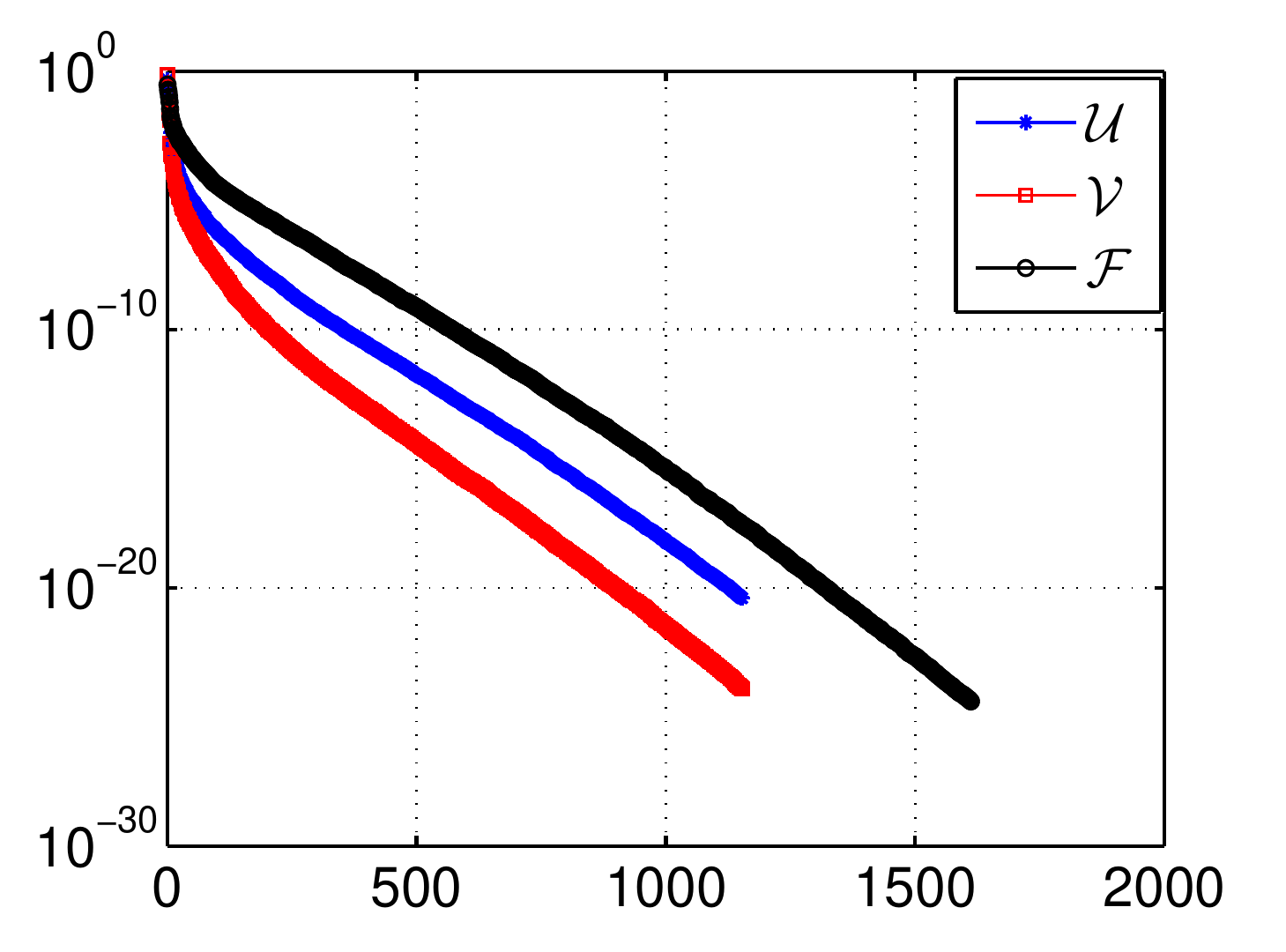}
\includegraphics[width=0.35\textwidth]{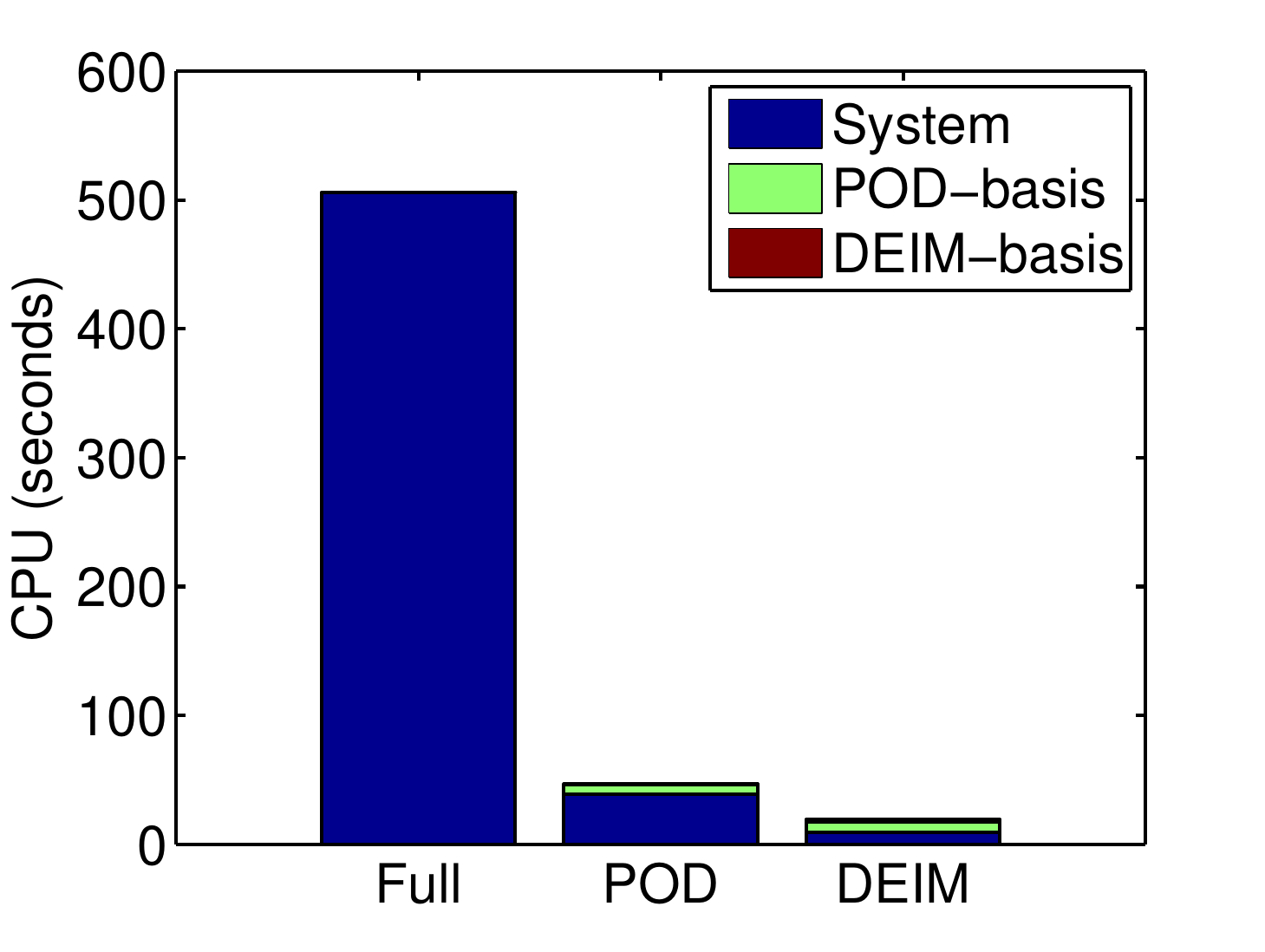}
\caption{(Left) Decay of the singular values of solution snapshots $\mathcal{U}$, $\mathcal{V}$ and of the nonlinear snapshots $\mathcal{F}$ ; (Right) CPU times for the computation of FOMs, POD and POD-DEIM ROMs for the parameter value $\mu=0.03$\label{cpu}.}
\end{figure}

\begin{figure}[htb!]
\centering
\includegraphics[width=0.31\textwidth]{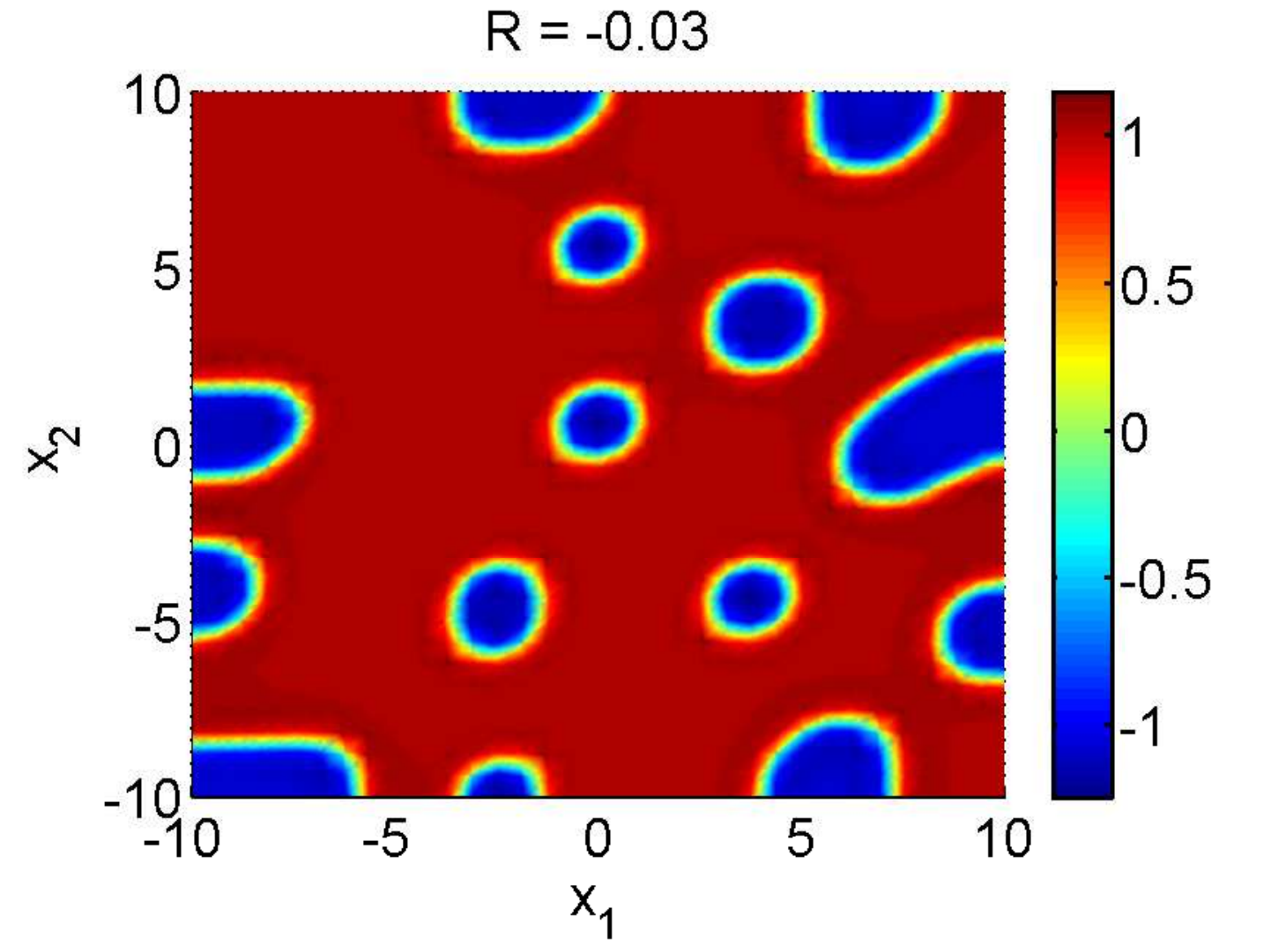}
\includegraphics[width=0.31\textwidth]{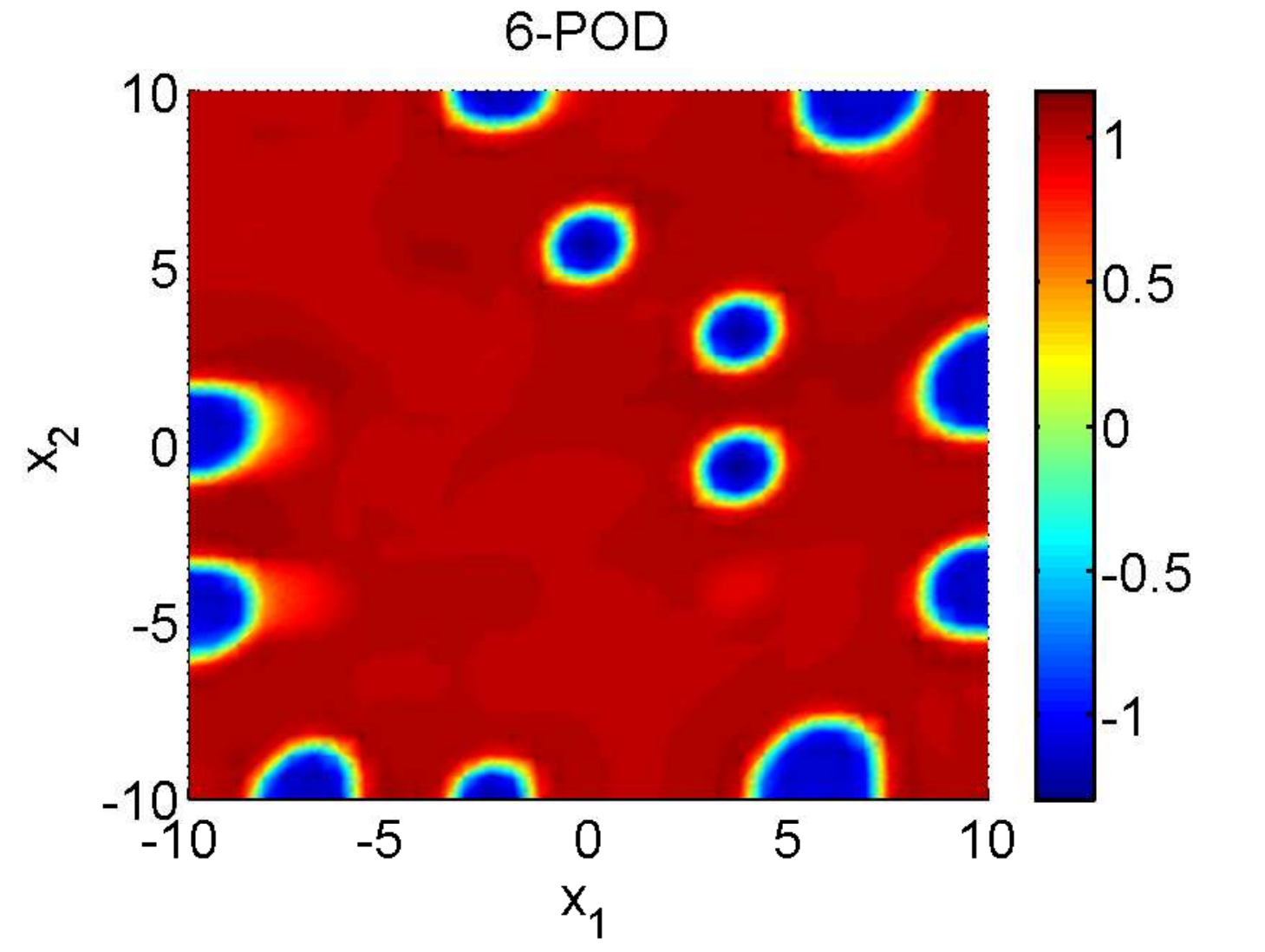}
\includegraphics[width=0.31\textwidth]{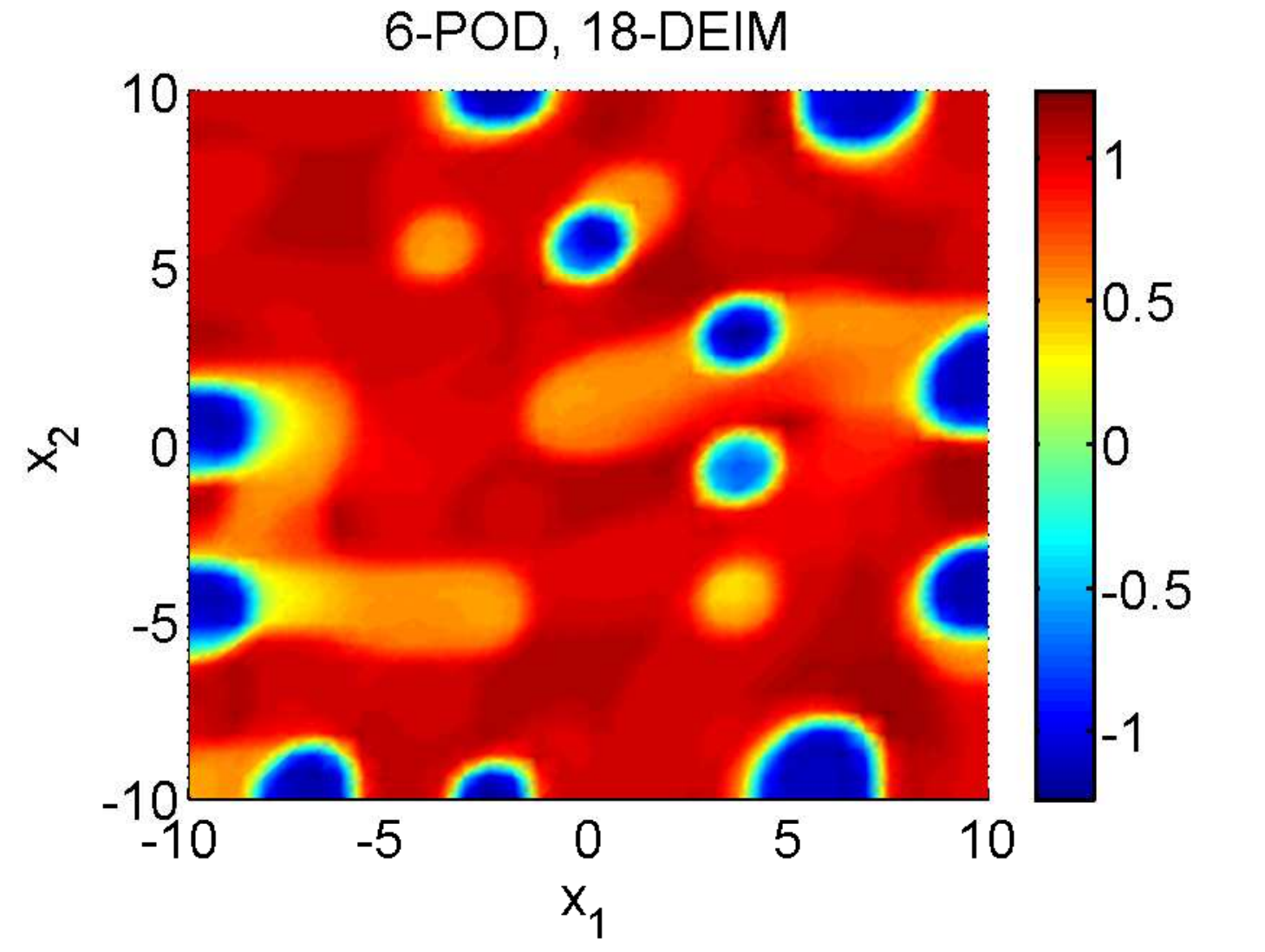}

\includegraphics[width=0.31\textwidth]{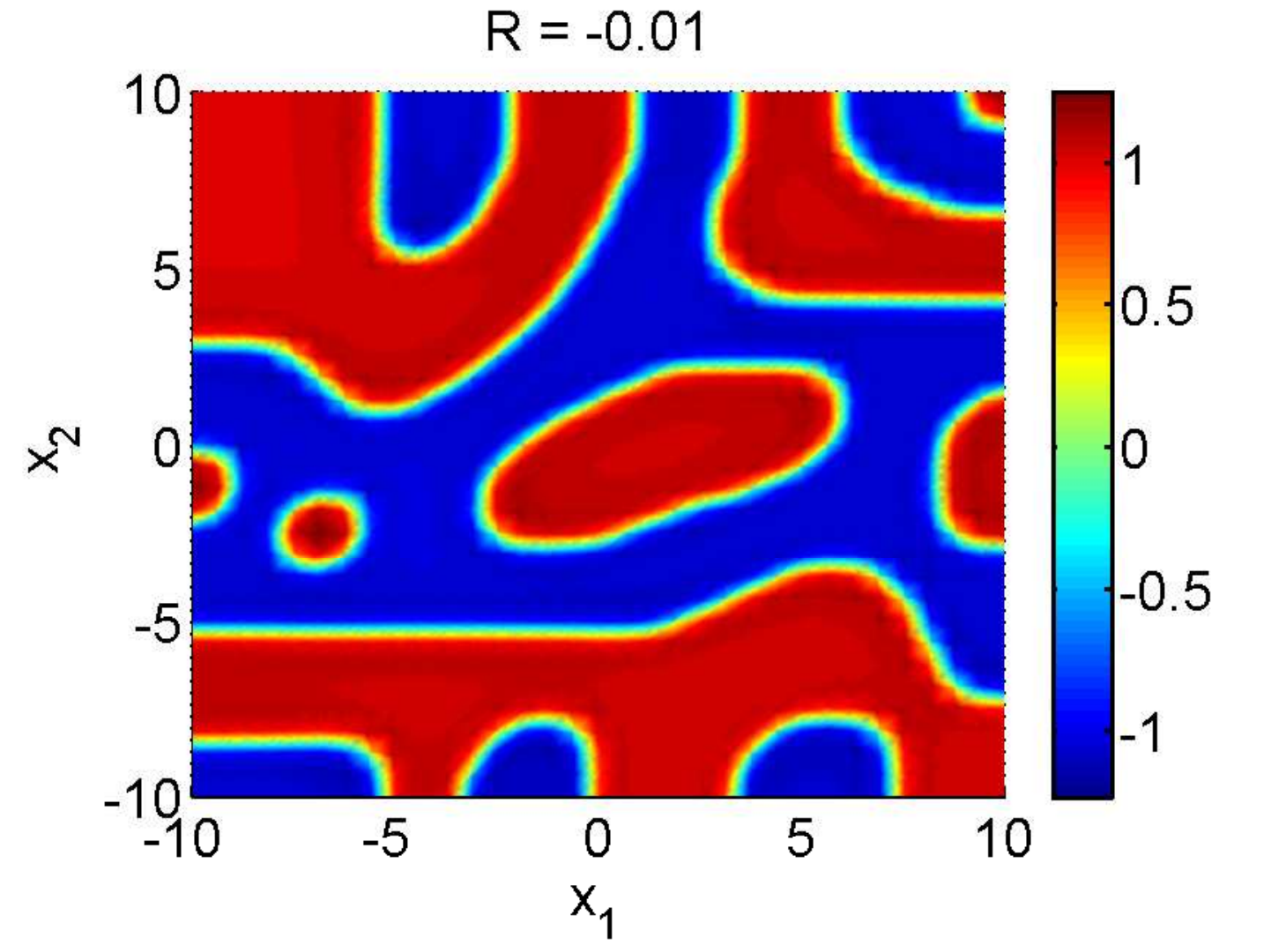}
\includegraphics[width=0.31\textwidth]{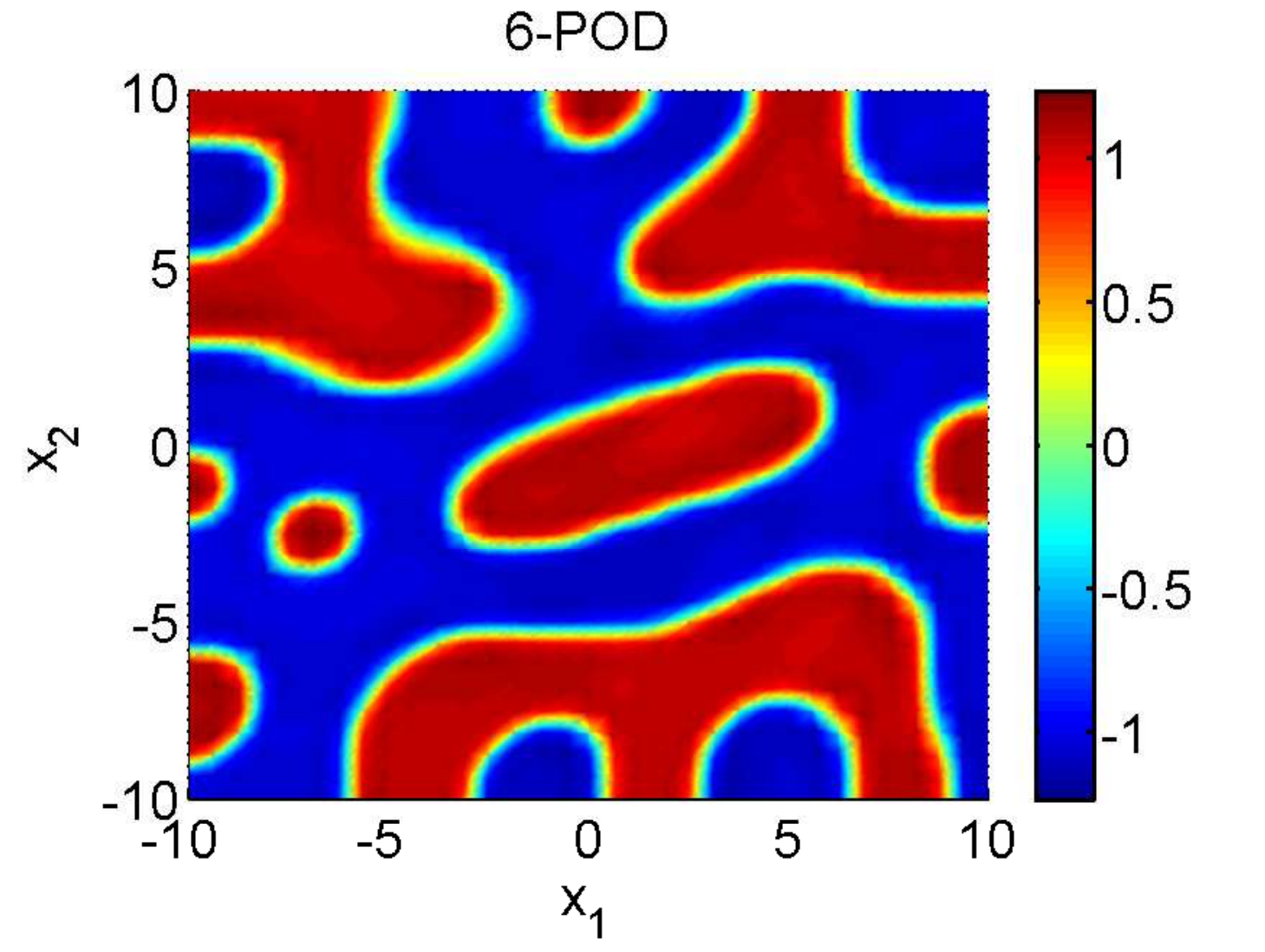}
\includegraphics[width=0.31\textwidth]{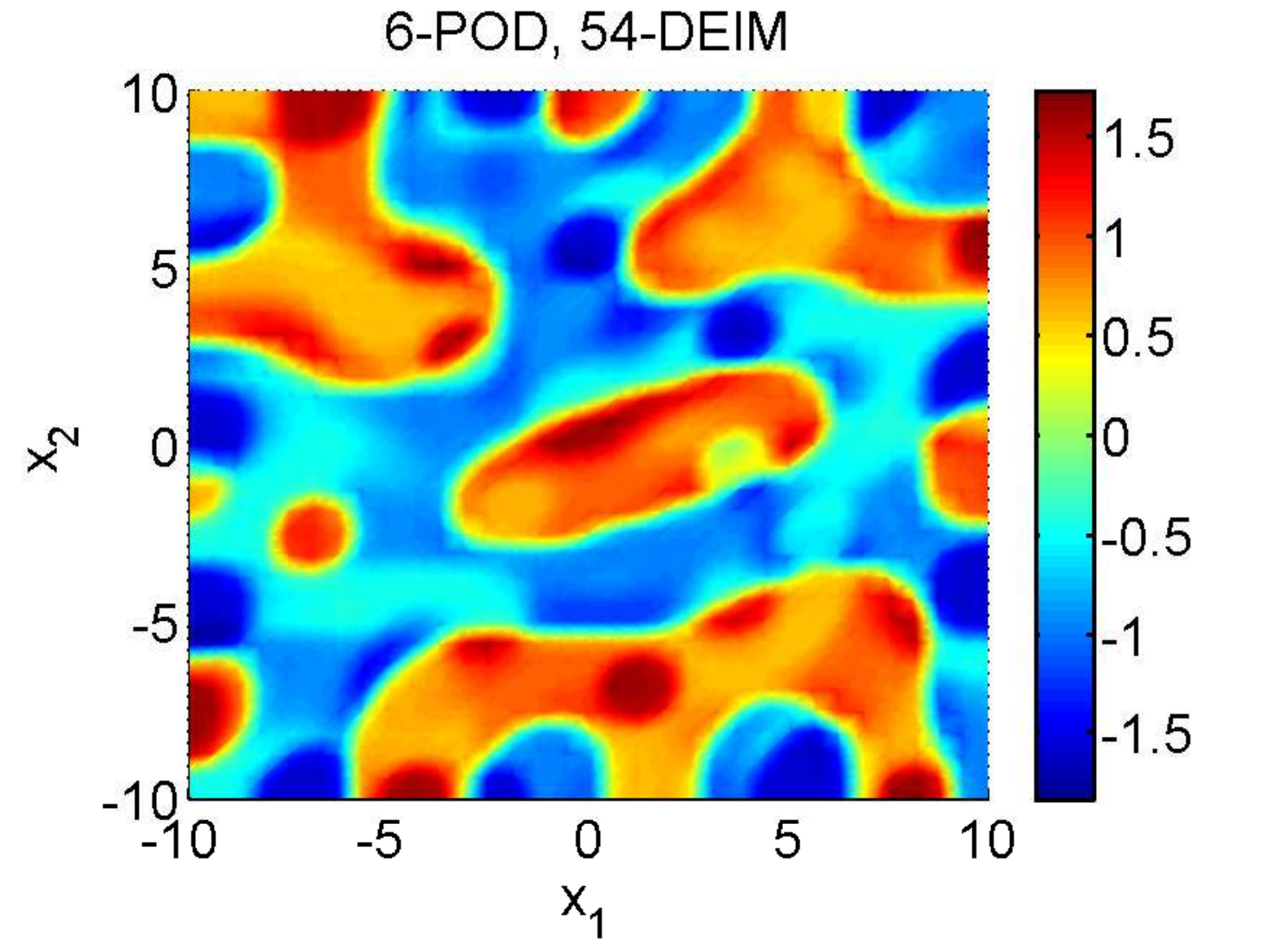}

\includegraphics[width=0.31\textwidth]{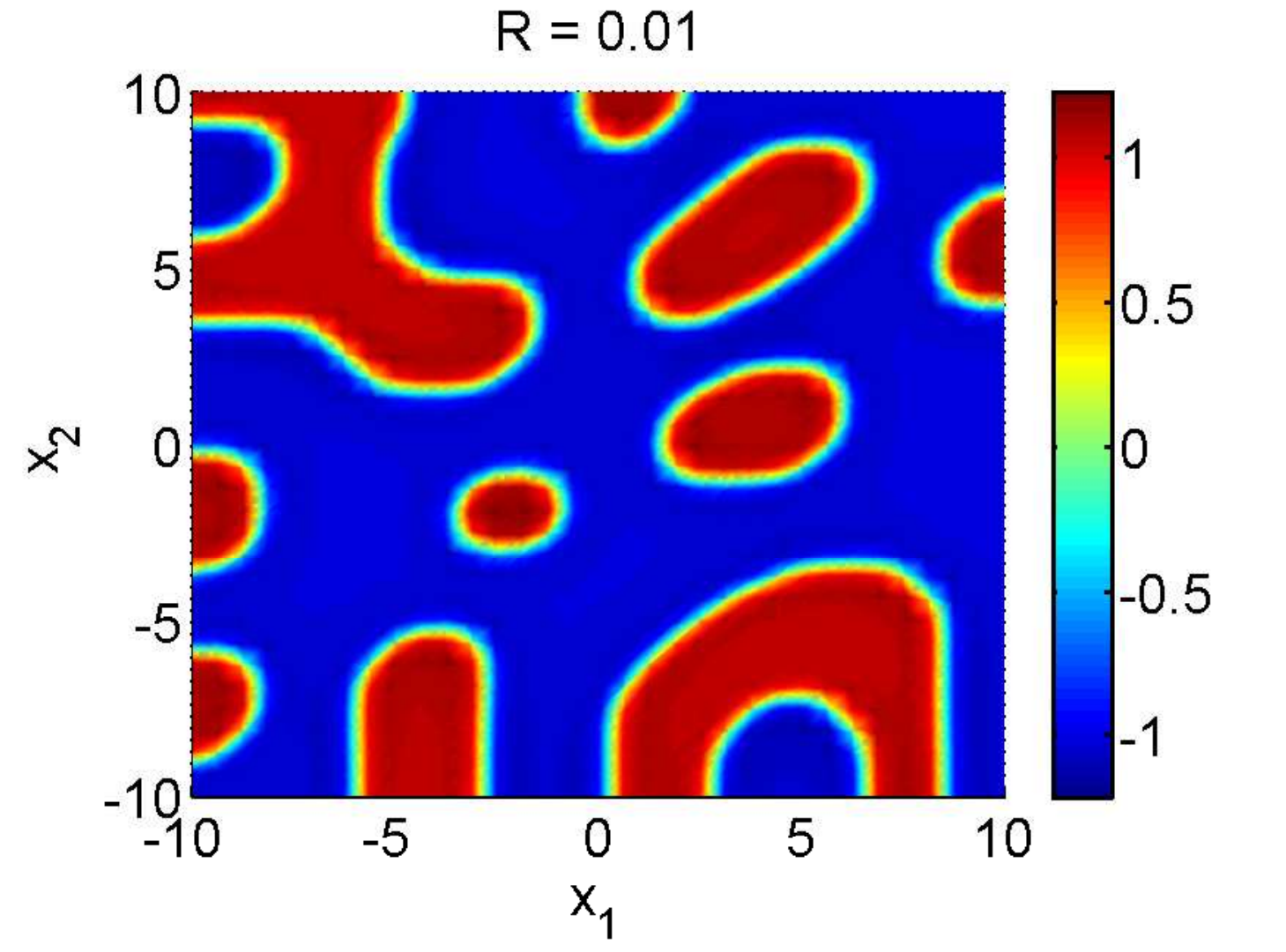}
\includegraphics[width=0.31\textwidth]{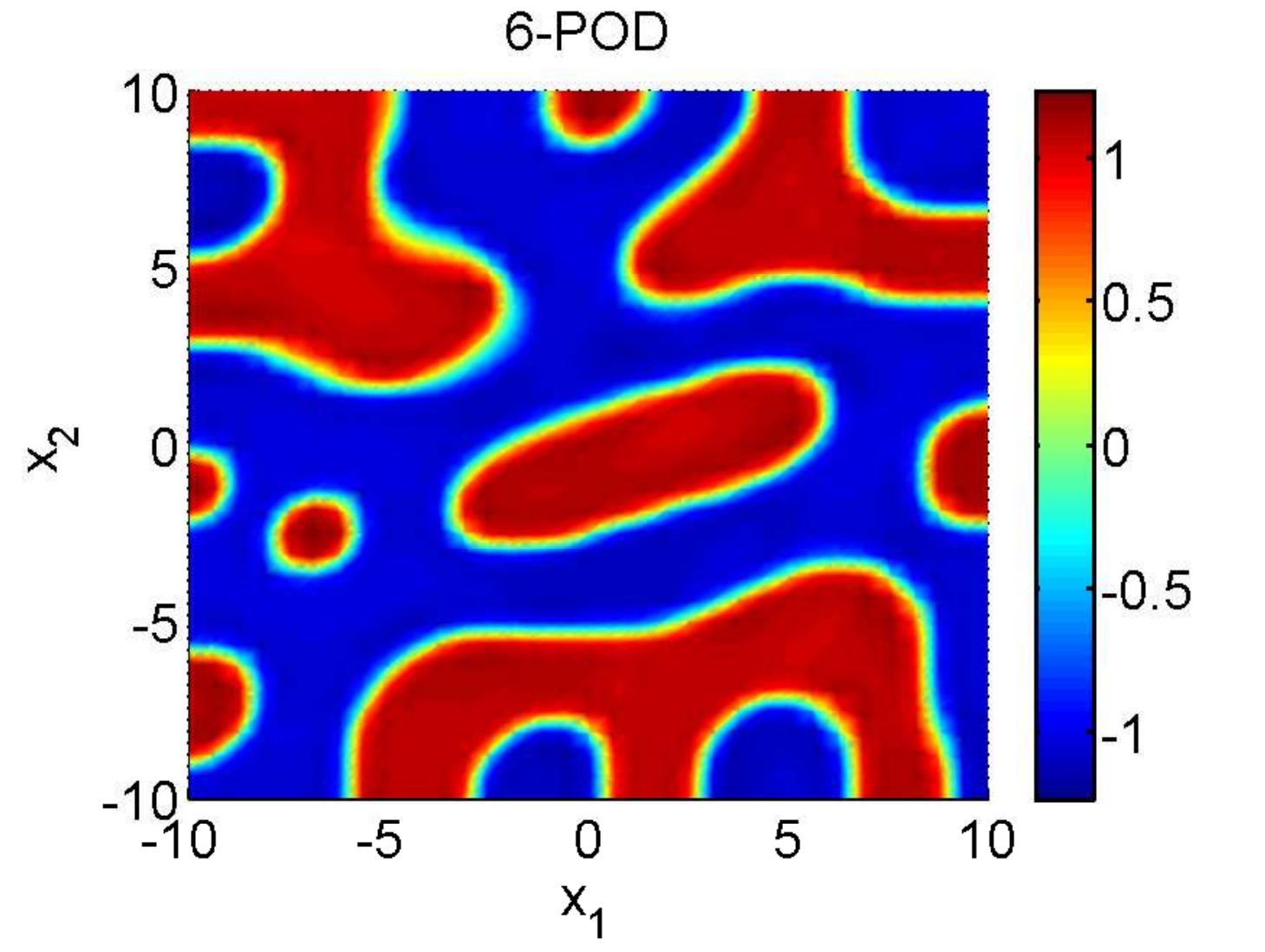}
\includegraphics[width=0.31\textwidth]{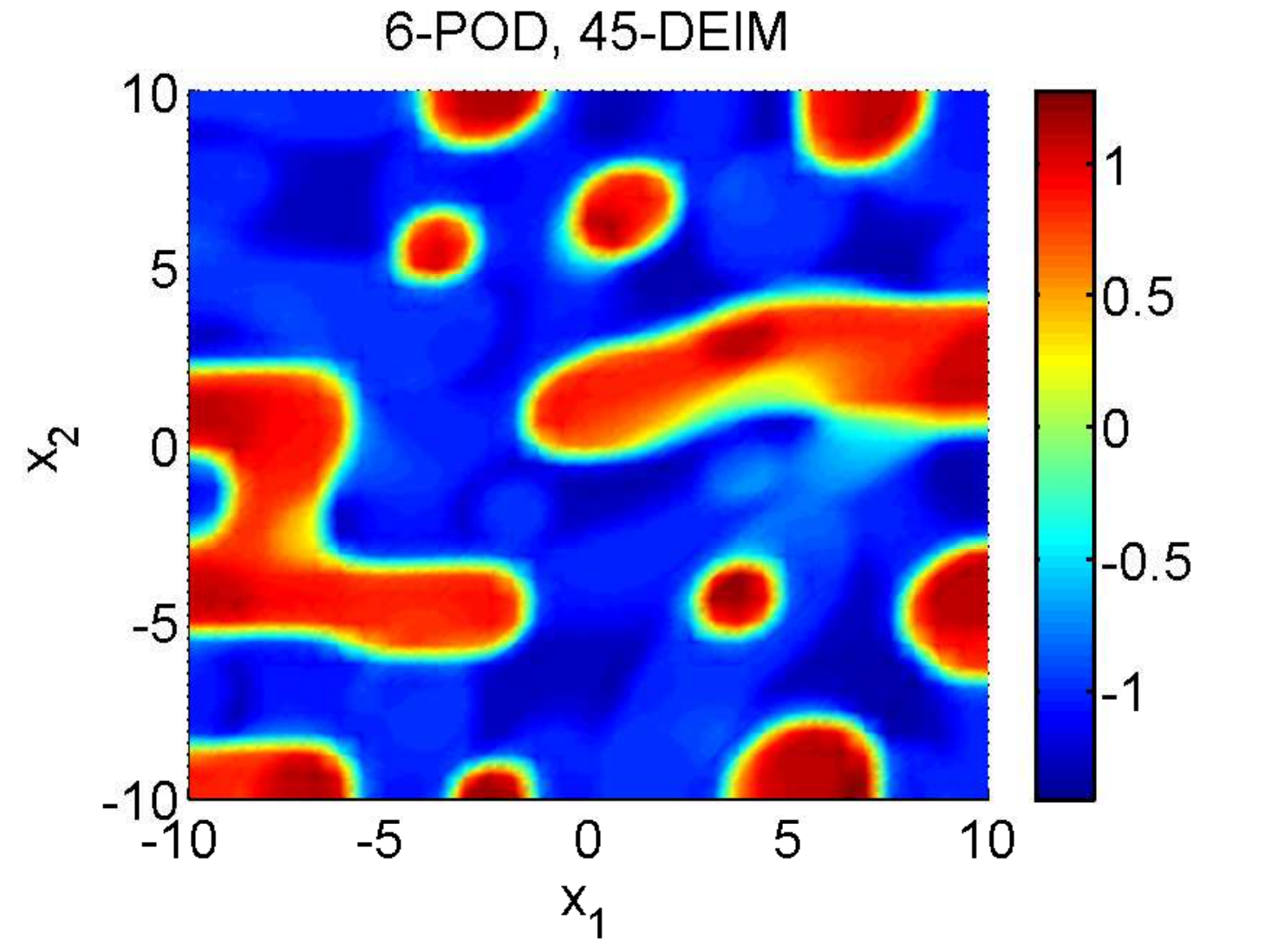}

\includegraphics[width=0.31\textwidth]{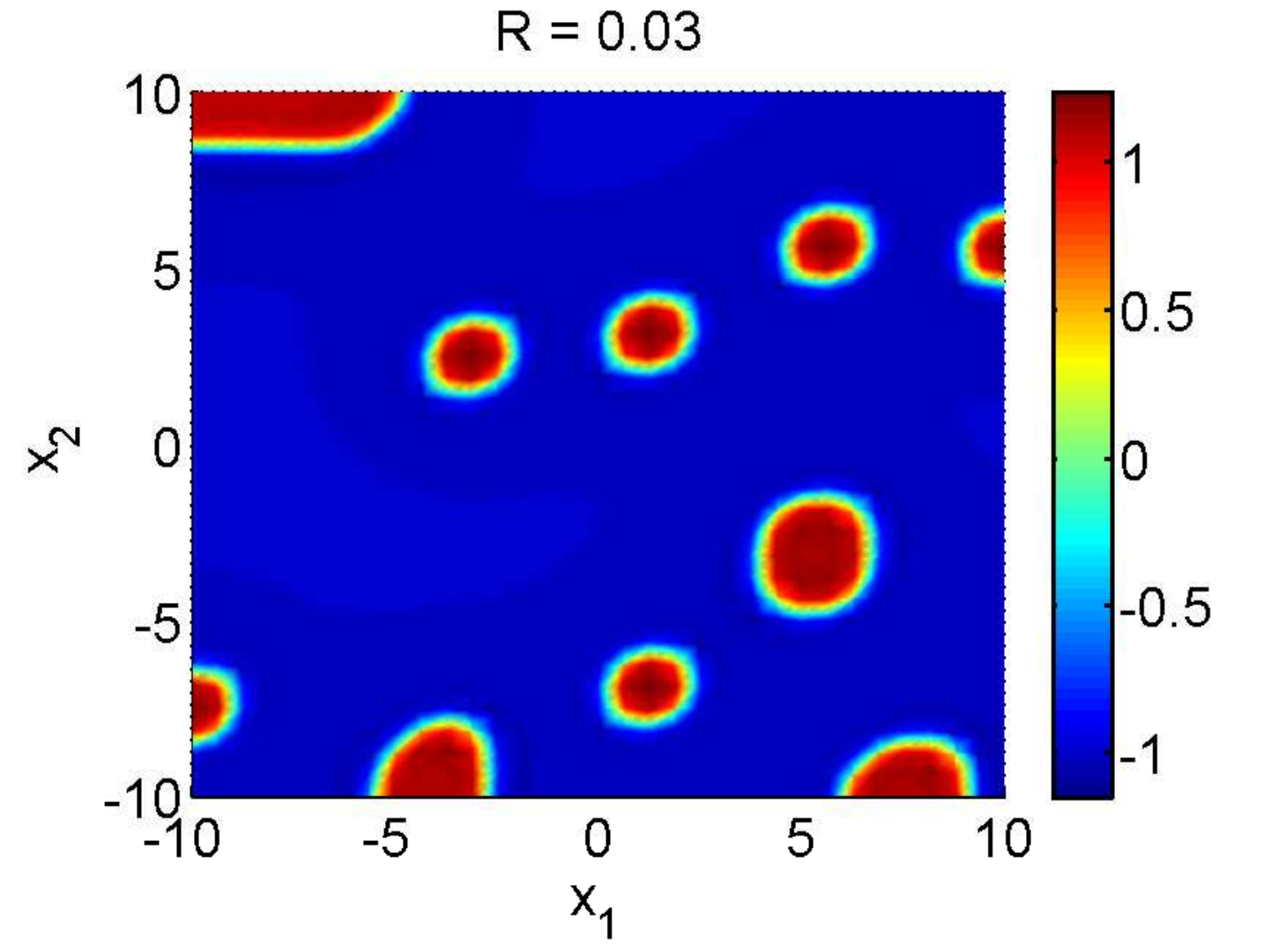}
\includegraphics[width=0.31\textwidth]{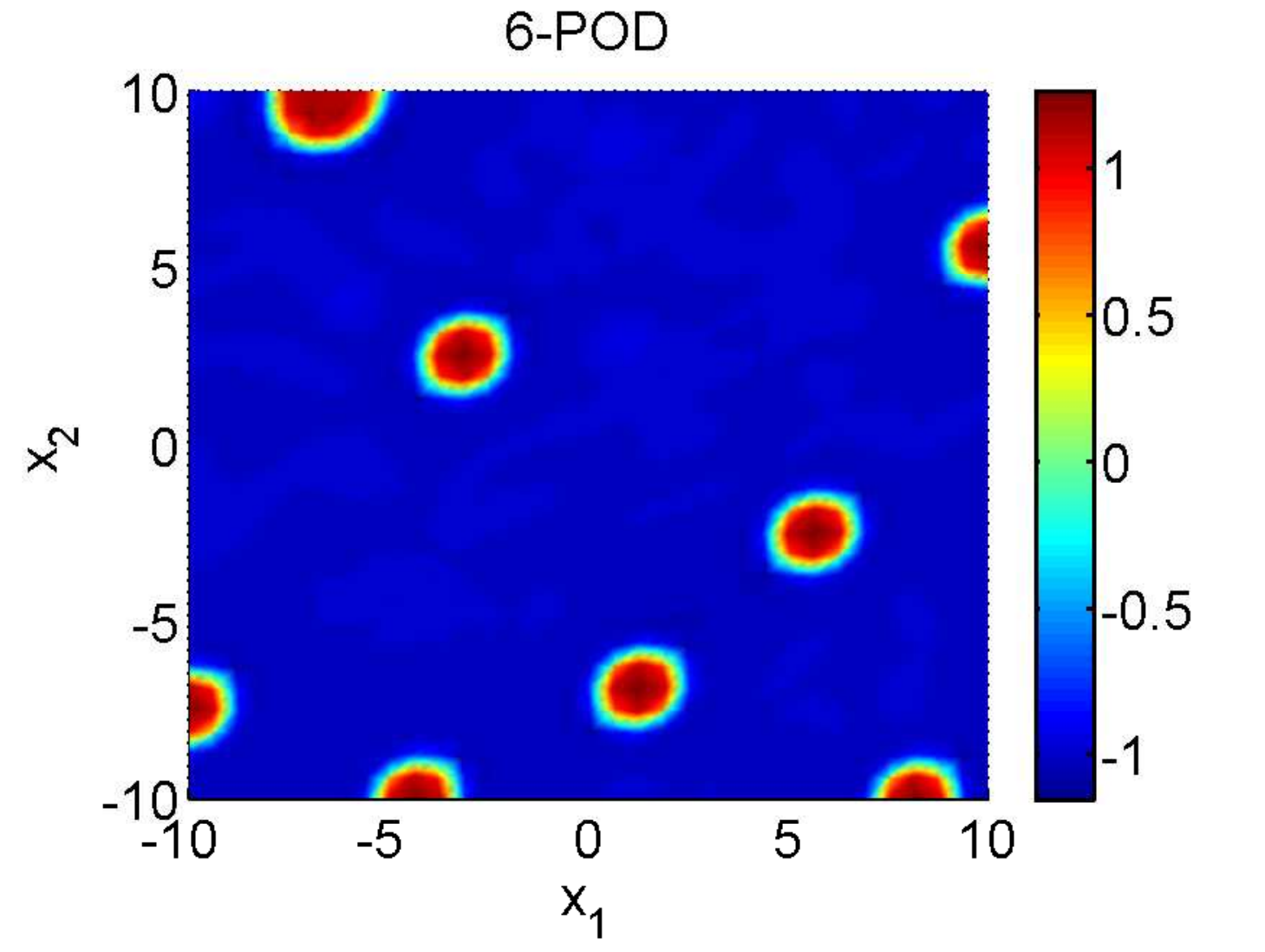}
\includegraphics[width=0.31\textwidth]{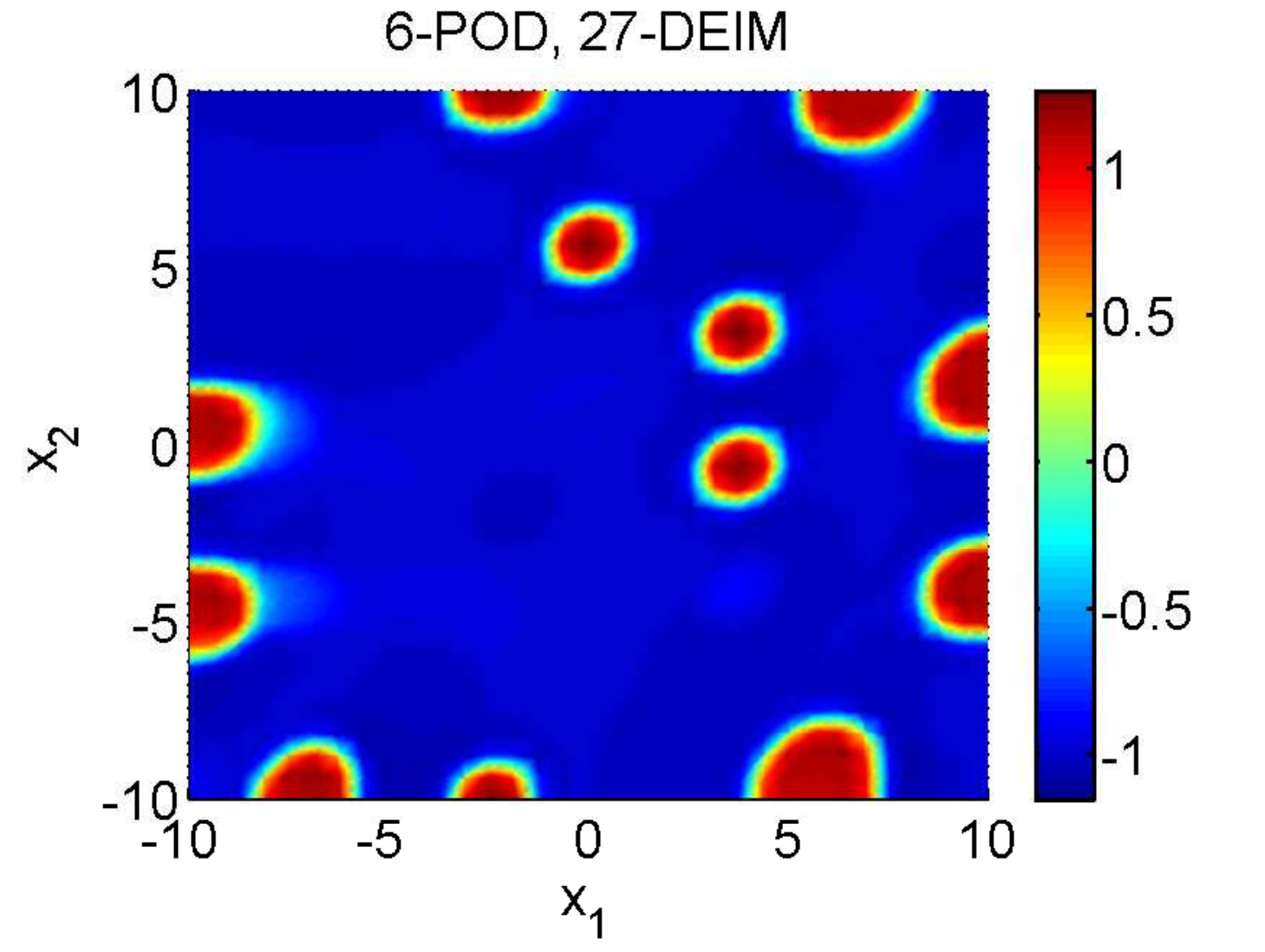}
\caption{Patterns for $u$ at the final time $T=1000$ with FOM (left), POD (middle) and POD-DEIM (right) for the parameter values $\mu\in\{-0.03,-0.01,0.01,0.03\}$ from top to bottom\label{pod}.}
\end{figure}

The decay of the singular values for the solution snapshots $\mathcal{U}$, $\mathcal{V}$ and nonlinear snapshots are given in Fig.~\ref{cpu}, left, and the CPU times of the FOMs and ROMs for the parameter value $\mu =0.03$ are shown in Fig.~\ref{cpu}, right. In Table~\ref{table1}  we give the CPU times  for FOMs, POD and POD-DEIM ROMs together with the speed-up factors $S_{POD}$ and $S_{DEIM}$, which demonstrate the efficiency of the DEIM. In Fig.~\ref{pod}, the patterns of FOMs, POD and POD-DEIM reduced solutions are shown at the final time $T=1000$. The ROM patterns in Fig.~\ref{pod} computed with POD are very close to those of the FOMs as in \cite{lago14}. But the patterns computed with POD-DEIM are less accurate than  those with the POD computed ones for some parameter values in  Fig.~\ref{pod}.  The DEIM does not improve the accuracy of the POD reduced model, but enormously reduces the computational complexity \cite{antil14}. The error bounds $\|(P^TW)^{-1}\|_2$ of moderate size for the DEIM approximations are also given in Table~\ref{table1}.

\begin{table}[htb!]
\centering
\caption{The computation times (in sec), speed-up factors $S_{POD}$ and $S_{DEIM}$, and the DEIM projection error bounds $\|(P^TW)^{-1}\|_2$}
\label{table1}
\begin{tabular}{ r | c |c |c| c |c |c }
$\mu$ & FOM & POD & POD-DEIM & $S_{POD}$ & $S_{DEIM}$ & $\|(P^TW)^{-1}\|_2$ \\
\hline
-0.03   &  527.3  & 34.5  & 7.5 & 15.31  & 70.21 &  28 \\
-0.01   &  501.9  & 33.4  & 13.2  & 15.05  &  38.08 &  33\\
0.01    &  522.3  & 32.9   & 11.9  & 15.88  & 43.67  &  41 \\
0.03    &  505.9  & 38.6   & 9.0   &  13.10  & 56.43  &  33 \\ \hline
\end{tabular}
\end{table}

\section*{Conclusions and Outlook}
We have demonstrated that the dG approximation can produce due to its local structure cost effective and accurate reduced order solutions by approximating the parameter dependent nonlinear terms with the DEIM. In a future work we will consider the parametrized FHN equation with the diffusivity coefficients $D_u$ and $D_v$ to compute the reduced order solutions by preserving the multiscale dynamics  of the activator $u$ and the inhibitor $v$  in time. Because the size of the SVD problem  can be prohibitive for the global POD, we will also apply the greedy POD.

\section*{Acknowledgments}
This work has been supported by METU BAP-07-05-2015 009.

\ifx\undefined\bysame
\newcommand{\bysame}{\leavevmode\hbox to3em{\hrulefill}\,}
\fi


\end{document}